\documentclass{elsart3-1}

\usepackage{amssymb}

\usepackage[english,francais]{babel}

\input xy
\xyoption{all}
\newtheorem{theorem}{Theorem}[section]
\newtheorem{theoremf}{Theorem}
\newtheorem{lemma}[theorem]{Lemma}
\newtheorem{e-proposition}[theorem]{Proposition}
\newtheorem{corollary}[theorem]{Corollary}
\newtheorem{corollaryf}[theoremf]{Corollary}
\newtheorem{e-definition}[theorem]{Definition\rm}

\renewcommand{\theequation}{\arabic{equation}}

\newcounter{sectionf}
\newcommand{\BR}{{\mathbb{R}}}

\def\bbr{{\Bbb R}}

\newcommand{\SL}{\mathrm{SL}}
\newcommand{\GL}{\mathrm{GL}}

\def\og{\leavevmode\raise.3ex\hbox{$\scriptscriptstyle\langle\!\langle$~}}
\def\fg{\leavevmode\raise.3ex\hbox{~$\!\scriptscriptstyle\,\rangle\!\rangle$}}

\journal{the Acad\'emie des sciences}
\begin{document}
\centerline{}
\begin{frontmatter}


\selectlanguage{english}
\title{Milnor-Wood inequalities for manifolds locally isometric to a product of hyperbolic planes}


\selectlanguage{english}
\author[authorlabel1]{Michelle Bucher},
\ead{mickar@math.kth.se}
\author[authorlabel2]{Tsachik Gelander}
\ead{tsachik.gelander@gmail.com}

\address[authorlabel1]{KTH Mathematics Department, 100 44 Stockholm, Sweden}
\address[authorlabel2]{Einstein Institute of Mathematics, The Hebrew University, Jerusalem, 91904, Israel}


\medskip

\begin{abstract}
\selectlanguage{english}
This note describes sharp Milnor--Wood inequalities for the Euler number of
flat oriented vector bundles over closed Riemannian manifolds
locally isometric to products of hyperbolic planes. One consequence is that
such manifolds do not admit an affine structure, confirming
Chern--Sullivan's conjecture in this case. The manifolds under consideration are of particular interest, since in contrary to many other locally symmetric spaces they do admit flat vector bundle of the corresponding dimension. 
When the manifold is irreducible and of higher rank, it is shown that flat oriented vector bundles are determined completely by the sign of the Euler number.

\vskip 0.5\baselineskip

\selectlanguage{francais}
\noindent{\bf R\'esum\'e} \vskip 0.5\baselineskip \noindent {\bf
In\'egalit\'es de Milnor--Wood pour vari\'et\'es localement
isom\'etriques \`a un produit de plans hyperboliques.} Nous
g\'en\'eralisons l'in\'egalit\'e classique de Milnor aux
vari\'et\'es localement isom\'etriques \`a un produit de plans
hyperboliques. Il en d\'ecoule que de telles vari\'et\'es
n'admettent pas de structure affine, confirmant dans ce cas la
conjecture de Chern--Sullivan. Contrairement \`a de nombreuses vari\'et\'es localement sym\'etriques, les vari\'et\'es consid\'er\'ees dans cette note admettent un fibr\'e vectoriel plat en dimension correspondante. Si les vari\'et\'es sont de plus
irr\'eductibles de rang sup\'erieur, nous montrons qu'un fibr\'e vectoriel orient\'e plat avec
nombre d'Euler non nul est, \`a
orientation pr\`es, unique. 

\end{abstract}
\end{frontmatter}

\selectlanguage{francais}
\section*{Version fran\c{c}aise abr\'eg\'ee}

Soit $\xi$ un $\mathrm{GL}^+(m,\bbr)$--fibr\'e principal, ou de fa\c
con \'equivalente, un fibr\'e vectoriel orient\'e de fibre
$\mathbb{R}^m$, sur une vari\'et\'e ferm\'ee $M$ orient\'ee de
dimension $m$. Rappelons que la classe d'Euler (r\'eelle) de $\xi$
est la classe de cohomologie $\varepsilon_m(\xi)\in
H^m(M,\mathbb{R})$ qui est l'image, par l'inclusion de coefficients
$\mathbb{Z}\hookrightarrow\mathbb{R}$, de l'obstruction \`a
l'existence d'une section non nulle dans le fibr\'e vectoriel
(associ\'e). Le nombre d'Euler de $\xi$ est le produit de Kronecker
de la classe d'Euler avec la classe fondamentale (r\'eelle) $[M]\in
H_m(M,\mathbb{R})$ de $M$:
$$\chi(\xi)=\langle \varepsilon_m(\xi),[M]\rangle.$$
Si $\xi$ admet une structure plate, c'est-\`a-dire si
$\xi$ est induit par une repr\'esentation $\rho:\pi_1(M)\rightarrow
\mathrm{GL}^+(m,\bbr)$ du groupe fondamental de $M$, alors il existe
une borne sur $|\chi(\xi)|$ ne d\'ependant que de $M$. En effet,
c'est une observation de Lusztig, que l'espace des repr\'esentations
de $\pi_1(M)$ dans $\mathrm{GL}^+(m,\bbr)$ est une vari\'e´t\'e
alg\'ebrique, et n'a en cons\'equence qu'un nombre fini de
composantes connexes. Comme le nombre d'Euler est constant sur les
composantes connexes, l'affirmation s'ensuit.

D\'enotons par $X$ l'espace hyperbolique r\'eel. Milnor, en premier,
a exhib\'e des fibr\'es plats dont la classe d'Euler est non nulle,
et de plus, donn\'e une borne optimale pour le nombre d'Euler de
fibr\'es plats au dessus de surfaces hyperboliques \cite{Mi58}. Nous
g\'en\'eralisons cette in\'egalit\'e aux vari\'et\'es localement
isom\'etriques \`a un produit de plans hyperboliques $X^n$, le cas
$n=1$ \'etant l'in\'egalit\'e de Milnor.

\begin{theoremf}\label{Theorem: French Milnor-Wood for products}
Soit $M$ une vari\'et\'e Riemannienne ferm\'ee localement
isom\'etrique \`a $X^n$ et soit $\xi$ un
$\GL^+(2n,\mathbb{R})$--fibr\'e sur $M$. Si $\xi$ admet une
structure plate, alors
\[
| \chi(\xi) |=| \left\langle \varepsilon_{2n}(\xi),[M]\right\rangle
| \leq\frac{1}{2^{n}}|\chi(M)|.
\]
\end{theoremf}

Il est imm\'ediat que si l'espace tangent $TM$ d'une vari\'et\'e
Riemannienne ferm\'ee $M$ localement isom\'etrique \`a $X^n$
admettait une structure plate, l'in\'egalit\'e
$|\chi(M)|=|\chi(TM)|\leq (1/2^n)|\chi(M)|$ d\'ecoulerait du
Th\'eor\`eme \ref{Theorem: French Milnor-Wood for products}, ce qui
n'est pas possible puisque $\chi(M)\neq 0$. Comme une structure
affine est clairement plate, nous en d\'eduisons une r\'eponse
affirmative partielle \`a la conjecture de Chern-Sullivan
pr\'edisant qu'une vari\'et\'e ferm\'ee avec charact\'eristique
d'Euler non nulle n'admet pas de structure affine:

\begin{corollaryf}\label{cor: french affine}
Une vari\'et\'e Riemannienne ferm\'ee localement isom\'etrique \`a
$X^n$ n'admet pas de structure affine.
\end{corollaryf}

Ce r\'esulat est nouveau pour $n>2$. Pour $n=1$ et $n=2$ cela
s'obtient d\'ej\`a par les in\'egalit\'es correspondantes dans
\cite{Mi58} et \cite{Bu08} respectivement. L'inexistence de
structure affine \textit{compl\`ete}, qui suivrait aussi de la
conjecture d'Auslander, est d\'emontr\'ee dans \cite{KoSu75}. Notons que m\^eme dans le cas d'un produit de surfaces,
l'inexistence de structure affine ne d\'ecoule pas directement de
l'inexistence de structure affine sur les surfaces. En effet,
Etienne Ghys nous a montr\'e un exemple d'une vari\'et\'e produit admettant une structure affine sans qu'aucun des facteurs n'en admettent une.

A rev\^etement fini pr\`es, les in\'egalit\'es du Th\'eor\`eme
\ref{Theorem: French Milnor-Wood for products} sont optimales,
c'est-\`a-dire que pour toute vari\'et\'e ferm\'ee $M$ localement
isom\'etrique \`a $X^n$, il existe un
$\GL^+(2n,\mathbb{R})$--fibr\'e $\xi$ admettant une structure plate
sur un rev\^etement $N$ fini de $M$ et $\chi(\xi)=(1/2^n)|\chi(N)|.$
Par contre, on ne peut pas r\'ealiser chaque entier dans
l'intervalle $[\frac{-|\chi(M)|}{2^n},\frac{|\chi(M)|}{2^n}]$ comme
nombre d'Euler de fibr\'e plat. En effet, nous pouvons raffiner le
th\'eor\`eme \ref{Theorem: French Milnor-Wood for products} comme
dans le Th\'eor\`eme \ref{thm: gnrl Euler number} ci-dessous. Le cas
extr\^eme \'etant quand $M$ est une vari\'et\'e que nous appellerons
rigide, c'est-\`a-dire que $M$ n'admet pas de rev\^etement fini se
d\'ecomposant en un produit contenant un facteur de dimension $2$.
Dans ce cas, les seules valeurs possibles pour le nombre d'Euler
d'un $\GL^+(2n,\mathbb{R})$--fibr\'e $\xi$ sur $M$ admettant une
structure plate sont $0$ et $\pm(1/2^n)|\chi(M)|$. De
plus, si le nombre d'Euler du fibr\'e plat est non nul, alors la
structure plate de $\xi$ est \`a orientation pr\`es unique. Voir
Th\'eor\`eme \ref{theorem: rigid} pour plus de d\'etails.

Observons enfin que contrairement au cas de dimension $2$, il n'est
pas possible de caract\'eriser les fibr\'es plats sur les
vari\'et\'es localement isom\'etriques \`a $X^n$ en fonction de leur
nombre d'Euler, car la classe d'Euler n'est pas un invariant complet
de classes d'isomorphies de $\GL^+(2n,\mathbb{R})$--fibr\'es en
dimension sup\'erieure. En effet, il n'est pas difficile de
construire deux fibr\'es avec m\^eme nombre d'Euler au dessus de
vari\'et\'es localement $X^n$, pour $n>1$, telle que l'un des
fibr\'e admette une structure plate, et l'autre non.

\selectlanguage{english}

\section{Historical introduction}
The first example of a nontrivial characteristic class of flat
bundle was given by Milnor in \cite{Mi58}, where he characterized
the $\GL^+(2,\mathbb{R})$--bundles admitting flat structures over
surfaces in terms of their Euler number: A
$\GL^+(2,\mathbb{R})$--bundle $\xi$ over a surface $\Sigma_g$ of
genus $g\geq 1$ admits a flat structure if and only if its Euler
number $\chi(\xi)$ satisfies the inequality
$|\chi(\xi)|=|\left\langle
\varepsilon_2(\xi),[\Sigma_g]\right\rangle|\leq g-1 $. Milnor's
inequality was later generalized to circle bundles by Wood
\cite{Wo75}.

In his groundbreaking essay \cite{Gr82}, Gromov naturally puts
Milnor's inequality in the context of bounded cohomology. Indeed,
canonical $L^1$ and $L^\infty$ norms can be defined on the spaces of
singular chains and cochains of a closed oriented $n$-dimensional
manifold $M$. These in turn induce seminorms on the respective real
valued homologies and cohomologies. It follows from Hahn-Banach theorem that
 \begin{equation}
|\left\langle \beta,[M]\right\rangle|=\left\Vert
\beta\right\Vert_\infty\cdot\left\Vert M\right\Vert,~\forall\beta\in
H^n(M)~\mathrm{with}~\left\Vert \beta \right\Vert_\infty<\infty,
\label{Equ: b(M)=product of norms}
\end{equation}
where $\left\Vert M\right\Vert$ denotes the $L^1$ seminorm of the
fundamental class of $M$, the so called {\it simplicial volume} of
$M$.
Thus, if $\beta$ is a characteristic class, a bound on the
characteristic number $|\left\langle \beta ,[M]\right\rangle| $ can
be obtained by bounding both $\left\Vert \beta\right\Vert _{\infty}$
and $\left\Vert M\right\Vert $. Unfortunately, estimating each of
these terms is usually very difficult.

Nonzero exact simplicial volume computations are rare. For oriented
surfaces $\Sigma_{g}$ of genus $g\geq1$, it is not difficult to show
that $\left\Vert \Sigma_{g}\right\Vert =2\left\vert
\chi(\Sigma_{g})\right\vert =4(g-1)$. In particular, if $g\ge 2$ and
$\Sigma_{g}$ is endowed with a hyperbolic structure, then
$\left\Vert \Sigma_{g}\right\Vert =\pi\cdot
\mathrm{Vol}(\Sigma_{g})$. More generally, if $M$ is an
$n$--dimensional closed hyperbolic manifold, then $\left\Vert
M\right\Vert =v_{n}\cdot \mathrm{Vol}(M)$ \cite{Gr82,Th78}, where
$v_{n}$ denotes the supremum of the volumes of geodesic simplices in
the $n$-dimensional hyperbolic space and is known explicitly in low
dimensions only. The only further computation of a nonzero
simplicial volume is given in \cite{Bu08} for manifolds locally
isomorphic to the product of two copies of the hyperbolic plane. In
this case, one has $\left\Vert M\right\Vert =6\cdot\chi
(M)=3/(2\pi^{2})\cdot \mathrm{Vol}(M)$.

It is known since \cite{Gr82} that characteristic classes of flat
$G$--bundles have finite $L^\infty$ seminorm when $G$ is a real
algebraic subgroup of $\mathrm{GL}(n,\mathbb{R})$, but actual upper
bounds for their norms are only known in special cases. For the
Euler class $\varepsilon_n$, Gromov \cite{Gr82} obtained from
Sullivan-Smillie's corresponding simplicial results that $\left\Vert
\varepsilon_{n}(\xi)\right\Vert _{\infty}\leq1/2^{n}$, whenever
$\xi$ is a $\GL^+(n,\mathbb{R})$--bundle admitting a flat structure.
Independently, Ivanov and Turaev \cite{IvTu82} exhibited an explicit
bounded cocycle representing the Euler class of flat bundles,
producing the same bound. In degree $2$, sharp upper bounds for the
K\"ahler class were computed by Domic and Toledo \cite{DoTo} in
terms of the rank of the associated symmetric space, later
generalized by Clerc and \O rsted \cite{ClOr} to include all
Hermitian symmetric spaces.

In view of the (im)possible seminorms computations, sharp
generalizations of Milnor's inequality were essentially carried
through in degree $2$ only. In this note, we announce some new sharp
upper bounds for the Euler number in degree $2n$ of every flat
$\GL^+(2n,\BR)$--bundle $\xi$ over a closed manifold $M$ which
admits a Riemannian structure locally isometric to a product of $n$
hyperbolic planes. Note however that our approach avoids estimating
the simplicial volume $\|M\|$ as well as the norm of the Euler
class. Instead, since by
(\ref{Equ: b(M)=product of norms}), $\left\Vert M\right\Vert=|\chi(M)|/\left\Vert \varepsilon_{2n}%
(TM)\right\Vert _{\infty}$, we prove Theorem \ref{Theorem: Milnor-Wood for products} by showing that
$\left\Vert \varepsilon_{2n}%
(\xi)\right\Vert _{\infty}\leq (1/2^n)\left\Vert \varepsilon_{2n}%
(TM)\right\Vert _{\infty}$ when $\xi$ is flat.


\section{Statement of results}

Denote by $X$ the real hyperbolic plane. We prove Milnor--Wood type
inequalities for manifolds admitting a Riemannian structure for
which the universal cover is isometric to the product of $n$
hyperbolic planes, in short, $X^n$--manifolds. The case $n=1$ is
Milnor's classical inequality.

\begin{theorem}\label{Theorem: Milnor-Wood for products}
Let $M$ be a closed $X^n$--manifold and $\xi$ a $\GL^+(2n,\mathbb{R})$--bundle over $M$. If $\xi$ admits a flat structure, then
\[
| \chi(\xi) |=| \left\langle \varepsilon_{2n}(\xi),[M]\right\rangle |
\leq\frac{1}{2^{n}}|\chi(M)|.
\]
\end{theorem}

It is an old conjecture of Chern (known also as the Chern--Sullivan conjecture), formulated independently by Milnor \cite{Mi58}, that a closed manifold of even dimension with nonzero
Euler characteristic cannot admit an affine structure. As a
consequence of Theorem \ref{Theorem: Milnor-Wood for products} we
derive the following partial affirmative answer:

\begin{corollary}\label{cor}
A closed $X^n$--manifold does not admit an affine structure.
\end{corollary}

In fact, we obtain the stronger statement that the tangent bundle
$TM$ of a closed $X^n$--manifold does not admit a flat structure.
Indeed, if it had, Theorem \ref{Theorem: Milnor-Wood for products}
would yield the impossible inequality 
$$
 | \chi(M) |=| \chi(TM) |=|
 \left\langle \varepsilon(TM),[M]\right\rangle |
 \leq(1/2^{n})|\chi(M)|.
$$  
Note that in general manifolds with nonzero Euler characteristic may admit a flat structure (cf. Smillie
\cite{Sm77}). 

The $X^n$--manifolds are of particular interest among all locally symmetric manifolds of even dimension. Indeed, while one can deduce directly from Margulis superrigidity theorem that many higher rank irreducible locally symmetric manifolds $M$ admit no flat bundle of dimension $\dim (M)$, $X^n$-manifolds do admit (in many cases a unique) such flat bundle (cf. Theorems \ref{thm: gnrl Euler number} and \ref{theorem: rigid}). By Theorem \ref{Theorem: Milnor-Wood for products} these bundles cannot be isomorphic to the tangent bundle $TM$.

Corollary \ref{cor} is new for $n> 2$. For $n=1$ and $n=2$ it
follows from the corresponding inequalities in \cite{Mi58} and
\cite{Bu08} respectively. The nonexistence of a \textit{complete}
affine structure is proved in \cite{KoSu75}. However, proving the
nonexistence of a {\it non}-complete affine structure is usually
much harder. Note also that even in the case of products of
hyperbolic surfaces, the nonexistence of affine structures on the
product can not be directly deduced from the nonexistence of affine structure
on the factors. Etienne Ghys showed us an example of a
product manifold which admits an affine structure,
while none of the factors does.

It is easy to construct examples of flat bundles over $X^n$--manifolds. More precisely, for any closed locally
$X^n$--manifold $M$, there exists a $\GL^+(2n,\mathbb{R})$--bundle
$\xi$ over a finite cover $N$ of $M$ such that $\xi$ admits a flat
structure and $\chi(\xi)=(1/2^n)|\chi(N)|.$ This shows that, up to finite cover, the inequality of Theorem \ref{Theorem:
Milnor-Wood for products} is sharp. However, not every
integer in the interval
$[\frac{-|\chi(M)|}{2^n},\frac{|\chi(M)|}{2^n}]$ is in general an
Euler number of a flat bundle, and Theorem \ref{Theorem: Milnor-Wood
for products} can be refined as in Theorem \ref{thm: gnrl Euler
number} below. We will say that an $X^n$--manifold is rigid if it
has no finite cover which decomposes as a product manifold with a
2-dimensional factor. This terminology is motivated by the (local,
Mostow and Margulis) rigidity theorems which apply for such
manifolds.

\begin{theorem} \label{thm: gnrl Euler number}
Let $M$ be a closed $X^n$--manifold, and $\xi$ a $\GL^+(2n,\mathbb{R})$-bundle over $M$. Let $N$ be a finite covering of $M$ of the form
$$N=\Sigma_{g_1}\times ...\times \Sigma_{g_k}\times N',$$
where $N'$ is rigid, $k\ge 0$, and the
$\Sigma_{g_i}$'s are surfaces of genus $g_i\geq 2$. Let $\xi_N$ be
the pullback of $\xi$ to $N$. If $\xi$ admits a flat structure, then
so does $\xi_N$ and
$$\chi(\xi_N)\in \{\frac{|\chi(N')|}{2^{\mathrm{Dim}(N')}}\Pi_{i=1}^k \ell_i :  |\ell_i|\leq g_i-1\}\cup\{0\}.$$
Moreover, upon passing to a finite cover of $N$, all the integers are actually attained as Euler numbers of flat bundles.
\end{theorem}

Note that contrarily to the $2$--dimensional case, one cannot
characterize flat bundles over closed $X^n$--manifolds by their
possible Euler numbers, since the Euler class is not a complete
invariant of isomorphism classes of $\GL^+(2n,\mathbb{R})$--bundles
in higher dimensions. Thus, the converses of Theorems \ref{thm:
gnrl Euler number} and \ref{Theorem: Milnor-Wood
for products} are in general not true. In contrast, for rigid manifolds the Euler number does characterize
flat bundles, whenever it does not vanish:
\begin{theorem}\label{theorem: rigid}
Let $M$ be a rigid closed $X^n$--manifold and let $\xi$ be a
$\GL^+(2n,\mathbb{R})$--bundle over $M$. Assume that $\xi$ admits a
flat structure. Then either $\chi(\xi)=0$ or
$|\chi(\xi)|=(1/2)^n|\chi(M)|$. Moreover, if $\chi(\xi)\ne 0$, then
$\chi(\xi)$ completely determines the bundle and the flat structure.
\end{theorem}


\section{About the proof of Theorem  \ref{Theorem: Milnor-Wood for products}}
Our proof of Theorem \ref{Theorem: Milnor-Wood for products}
combines bounded cohomology, representation theory and rigidity theory.

Let $\xi$ be a $\GL^+(2n,\mathbb{R})$--bundle with nonzero Euler number
admitting a flat structure over a closed $X^n$--manifold $M$. Upon
passing to a double cover of $M$ we can suppose that $M$ is
oriented. The flat bundle $\xi$ is induced by a representation
$$\rho:\pi_1(M)\longrightarrow\GL^+(2n,\mathbb{R}).$$
Denote by $\varepsilon_{2n}\in H^{2n}(\prod_{i=1}^n\GL(2,\mathbb{R})^+)$ the (universal) Euler class of flat bundles, and note that by naturality, the Euler class of $\xi$ is the pullback  $\varepsilon_{2n}(\xi)=\rho^*(\varepsilon_{2n})$ of $\varepsilon_{2n}$ by $\rho$. The main part of our proof consists of several steps which allow us
to reduce the general case, without changing the $L^\infty$ seminorm
of the Euler class $\rho^*(\varepsilon_{2n})$, to the situation
where $\mathrm{Im}(\rho)$ is contained in a direct product of $n$ diagonally embedded $GL^+(2,\mathbb{R})$'s. We prove and make use of the following general lemma:

 \begin{lemma}\label{lem:prod}
Let $G_{i}$ be groups and let
$\rho:\prod_{i=1}^{n}G_{i}\to
\mathrm{GL}_{m}^{+}({\mathbb{R}})$ be an orientable representation
such that $\rho(G_{i})$ is not amenable for each $i$. Then
$m\geq2n$, and in case $m=2n$, up to replacing each $G_{i}$ by
a finite index subgroup, $V={\mathbb{R}}^{2n}$ decomposes as a
direct sum of $n$ invariant planes $V_{i}$ such that
$\rho(G_{i})$ acts irreducibly on $V_{i}$ and scalarely (i.e.
by multiplications by scalar diagonal matrices) on each $V_{j}, ~j
\neq i$. 
\end{lemma}

Upon replacing $M$ by a finite cover, we can without loss of generality suppose that we are in one of the following three cases:
\begin{enumerate}
\item $M$ is completely reducible, i.e. a finite cover of $M$ is a product of $n$ surfaces $\Sigma_{g_1}\times ... \times \Sigma_{g_n}$. In this case, we applying Lemma \ref{lem:prod} directly to the product $\prod_{i=1}^n\pi_1(\Sigma_{g_i})$ and conclude that if $\left\Vert \rho^*(\varepsilon_{2n})\right\Vert \neq 0$, then none of the $\rho(\pi_1(\Sigma_{g_i}))$'s is amenable and the representation factors, up to conjugation and upon replacing $\pi_1(M)$ by a finite index subgroup, through the diagonally embedded product  $\Pi_{i=1}^n\GL^+(2,\mathbb{R})$. 
\item $M$ is rigid, i.e. no finite cover of $M$ admits a two dimensional factor. Here, we use Margulis Superrigidity Theorem to show that $\rho(\pi_1(M))$ is, up to finite index, contained in a connected semisimple Lie group $S$ for which all the noncompact simple factors are locally isomorphic to $\mathrm{PSL}(2,\mathbb{R})$. Furthermore, if $S$ has $n$ factors which are all noncompact, then we get from Lemma \ref{lem:prod} applied to $S$ that $S$ is conjugated to  the diagonally embedded product  $\Pi_{i=1}^n\mathrm{SL}(2,\mathbb{R})$. 
\item $M$ is mixed, i.e. admits a finite cover which is a direct product of surfaces and a rigid manifold. In this last case, the image $\rho(\pi_1(M))$ is contained in a semidirect product $S\ltimes A$, where $S$ is semisimple and $A$ is a connected amenable normal subgroup. Since $A$ is amenable, denoting by $p$ the projection $p:S\ltimes A \rightarrow S$, we get $\left\Vert \rho^*(\varepsilon_{2n})\right\Vert_\infty =\left\Vert \rho^*(p^*(\varepsilon_{2n}))\right\Vert_\infty $, and we can hence replace the representation $\rho$ by the composition of $\rho$ with the projection $p$. If for one of the factors $M_i$ of $M$, the image $\rho(\pi_1(M_i))$ lies in a closed amenable group, or in the group of isometries of a symmetric space of lower dimension than the dimension of $M_i$, then $\rho^*(\varepsilon_{2n})=0$. Otherwise, it again follows from Lemma \ref{lem:prod} that $S$ has to be conjugated to a subgroup of the diagonally embedded product  $\Pi_{i=1}^n\SL(2,\mathbb{R})$. 
\end{enumerate}

To summarize, we have reduced the proof to the case where the representation $\rho$
factors through
\[
\xymatrix{\pi_1(M)\ar[r]^{\rho \ \ \ \ \ }\ar[dr]_{\rho_0} & \GL^+(2n,\bbr) \\
                 & \Pi_{i=1}^n\GL^+(2,\mathbb{R}).\ar@{^{(}->}[u]^i      }
\]
In view of Whitney's product formula for the Euler class, the pullback $i^*(\varepsilon_{2n})\in H^{2n}(\prod_{i=1}^n \GL(2,\mathbb{R})^+)$ of the Euler class is the cup product $\varepsilon_{2}\cup
... \cup \varepsilon_{2}$ of the Euler class on the factors and hence
\begin{equation}
\left\Vert\rho^*(\varepsilon_{2n})\right\Vert_\infty=\left\Vert\rho_0^*(i^*(\varepsilon_{2n}))\right\Vert_\infty=
\left\Vert\rho_0^*(\varepsilon_{2}\cup ... \cup
\varepsilon_{2})\right\Vert_\infty\leq\left\Vert\varepsilon_{2}\cup
... \cup \varepsilon_{2}\right\Vert_\infty.\label{Equ: rho(eps) <
eps U ... U eps}
\end{equation}
For the simplicial volume of $M$, we show using Hirzebruch's Proportionality Principle \cite{Hi58} that
\begin{equation}
\left\Vert M\right\Vert =\frac{|\chi(M)|}{2^n\left\Vert
\varepsilon_{2}\cup ... \cup
\varepsilon_{2}\right\Vert_\infty}.\label{Equ: simpl vol=X(M)/2^n
Euler norm}
\end{equation}
Finally, combining (\ref{Equ: rho(eps)
< eps U ... U eps}) and (\ref{Equ: simpl vol=X(M)/2^n Euler norm}), and using the duality of the $L^1$ and $L^\infty$ seminorms
(\ref{Equ: b(M)=product of norms}),
we obtain
$$
|\chi(\xi)|=\left\Vert\rho^*(\varepsilon_{2n})\right\Vert_\infty\left\Vert
M\right\Vert\leq \frac{1}{2^n}|\chi(M)|.
$$
The details will appear elsewhere.





\end{document}